\documentclass[conference,harvard,brazil,english,portuguese]{sbatex}

\makeatletter
\def\verbatim@font{\normalfont\ttfamily\footnotesize}
\makeatother

\usepackage[utf8]{inputenc}

\usepackage{icomma} 
\usepackage{amsmath,amssymb,amsfonts} 
\usepackage{verbatim}

\usepackage{graphicx} 
\usepackage{psfrag} 
\usepackage{float} 
\usepackage{multirow}

\usepackage{url} 

\begin{document}

\title{SIMULAÇÃO DE SISTEMAS DINÂMICOS COM ANÁLISE INTERVALAR: UM ESTUDO DE CASO COM O CIRCUITO RLC}

\author{Peixoto, M. L. C.}{marciapeixoto93@hotmail.com}
\address{GCOM - Grupo de Controle e Modelagem \\ Departamento de Engenharia Elétrica \\ UFSJ - Universidade Federal de São João del-Rei\\ Pça. Frei Orlando, 170 - Centro - 36307-352 - São João del-Rei, MG, Brasil}

\author[1]{Nepomuceno, E. G.}{nepomuceno@ufsj.edu.br}

\author[1]{Junior Rodrigues, H. M.}{heitormrjunior@hotmail.com}

\author[1]{Martins, S. A. M.}{martins@ufsj.edu.br}

\author[1]{ Amaral,G. F. V.}{amaral@ufsj.edu.br}

\twocolumn[

\maketitle


\selectlanguage{english}
\begin{abstract}

Differences between computer simulation of dynamical systems and laboratory experiments are common in teaching and research in engineering. Normally,  numerical inaccuracy and the non-ideal behaviour of the devices involved in the experiment are the most common explanations. With the application of interval analysis, it is possible to incorporate the numerical and parametric uncertainties in the simulation, allowing a better understanding of the play between simulation and experiment. This article presents a case study in which an step input is applied to an RLC circuit. Using the toolbox Intlab for Matlab, it was possible to present a computer simulation with the range that encompasses the experimental results . Comparison of simulation with experimental data show the success of the technique and indicates a potential content to be delivered to undergraduate engineering courses.

\end{abstract}

\keywords{Dynamical Systems, 
Propagation of Errors, Interval analysis, Intlab toolbox.}

\selectlanguage{brazil}
\begin{abstract}

Diferenças entre a simulação computacional de sistemas dinâmicos e
experimentos em laboratório são comuns no ensino e pesquisa em Engenharia. Normalmente, incertezas numéricas da simulação e o comportamento não-ideal dos componentes envolvidos no experimento são as explicações mais comuns.
Com a aplicação da análise intervalar, é possível incorporar as incertezas numéricas e paramétricas na simulação, permitindo uma melhor compreensão entre os limites da simulação com os dados coletados do experimento. Este artigo apresenta um estudo de caso em que uma entrada ao degrau é aplicada a um circuito RLC. Com o uso do toolbox Intlab para Matlab, foi possível apresentar uma simulação computacional com o intervalo que englobe os resultados experimentais. A comparação da simulação com dados experimentais, mostrou o êxito da técnica e indica um potencial conteúdo a ser ministrado para cursos de graduação em Engenharia.

\end{abstract}

\keywords{Sistemas Dinâmicos, Propagação de erros, Análise Intervalar, toolbox Intlab.}
]

\selectlanguage{brazil}

\section{Introdução}

A aprendizagem  por meio de experimentos laboratoriais é uma parte fundamental
do ensino de engenharia \cite{Feisel2005,Wollenberg2010}. Neste contexto, a realização de experimentos que caracterizam o comportamento de sistemas são bem comuns no ensino e educação dos engenheiros e especialistas na área \cite{Murray2002}. Comumente  estes experimentos são realizados com uma entrada em degrau, pois a partir desta é possível  estimar os parâmetros do sistema,  bem como projetar controladores.

No entanto, em Engenharia, incertezas aparecem frequentemente em medidas e modelagem. Por mais preciso e avançado que seja um determinado dispositivo de medição, geralmente haverá um erro associado. Este erro é normalmente conhecido como incerteza da medição.  Na modelagem de um fenômeno físico, existem incertezas relacionadas à diferença entre a realidade e o modelo escolhido, à impossibilidade de quantificar perfeitamente os parâmetros de um modelo e em entender o fenômeno \cite{Kulisch2012}. Existem técnicas para reduzir esta incerteza, como por exemplo repetir a medição várias vezes e calcular a média aritmética dos valores medidos \cite{Bich2006}.

No âmbito da Engenharia Elétrica e cursos afins, as aproximações estão presentes em componentes elétricos, que apresentam um valor nominal que corresponde ao provável valor e uma tolerância que indica uma possível variação do valor do componente em torno de seu valor nominal.  
Por exemplo, considerando um resistor com valor nominal de 100 $\Omega$ com tolerância de 5\%, sabe-se que seu valor real pode estar entre 95 e 105 $\Omega$.
Da mesma forma acontece com os capacitores. Por exemplo, o capacitor cerâmico 103M  apresenta um valor nominal de capacitância de 10 nF com uma tolerância de 20\%. 

\begin{figure}[htb]
	\centering
	\includegraphics[width=7.5cm,height=6cm]{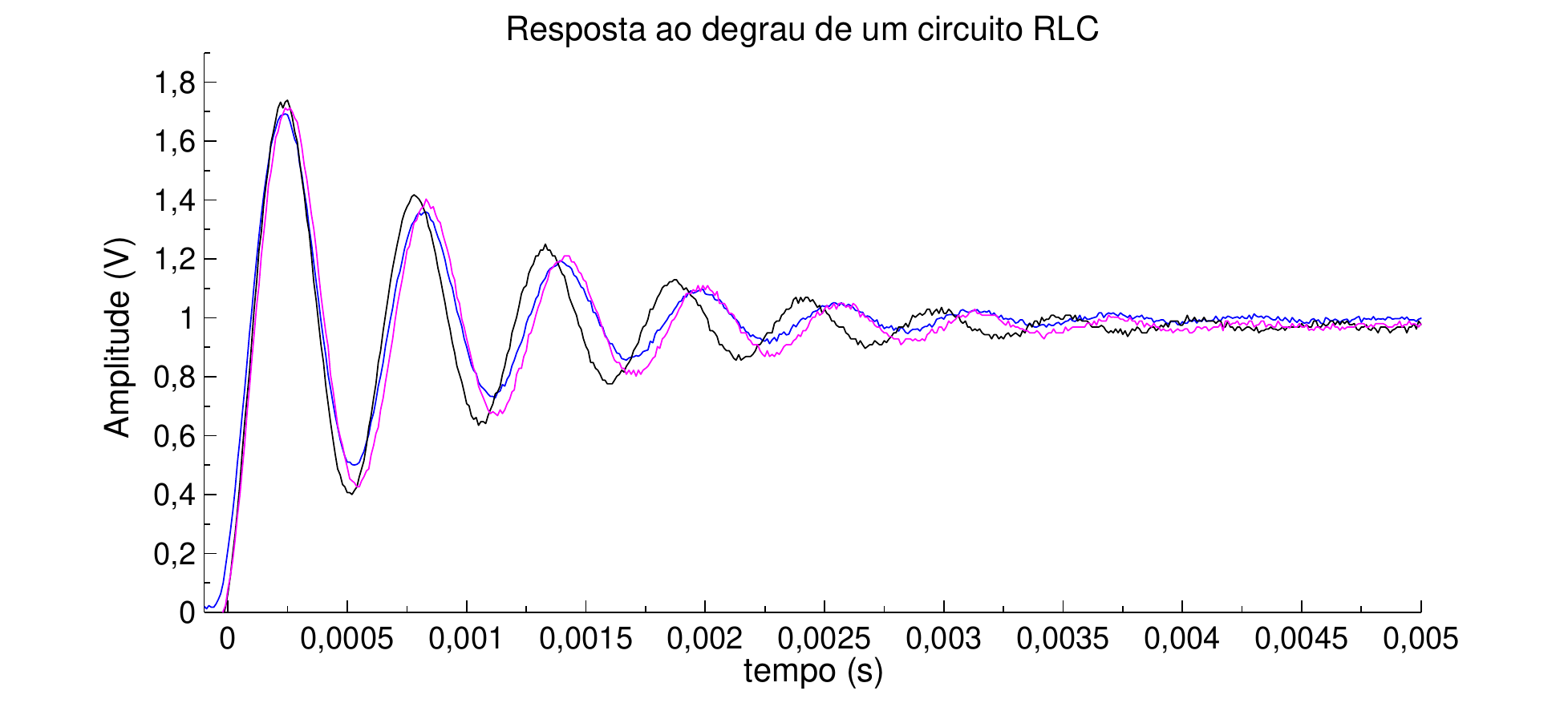}
	\caption{Exemplo de incertezas na resposta ao degrau de um circuito RLC. Três experimentos, em momentos diferentes, do mesmo circuito mostram resultados ligeiramente diferentes.}
	\label{fig:prat}
\end{figure}

A partir da  Figura \ref{fig:prat}, é possível perceber como as variações afetam um sistema. Nesta figura há três medições para a entrada ao degrau a um circuito RLC série, realizadas em instantes diferentes. Como os valores dos componentes não podem ser determinados exatamente, devido à variação de temperatura, umidade, precisão dos instrumentos, dentre outros fatores, é possível notar que o comportamento geral das três curvas é similar, entretanto não coincidem totalmente. A intensidade desta variação na saída pode ser considerável em sistemas mais complexos.
O tratamento adequado dos erros de medição é um aspecto crucial do procedimento laboratorial cuidadoso. Um estudante deve ter as ferramentas para identificar as fontes dos erros experimentais, estabelecer incertezas nas medições e propagar essas incertezas através dos cálculos que levam aos resultados finais \cite{Rothwell2012}.

Para melhor compreensão dos temas abordados nos laboratórios, as simulações computacionais têm sido muito utilizadas para diversas tarefas de ensino. As simulações em um ambiente de ensino trazem diversos benefícios pois  existe o contato direto com o modelo matemático do sistema dinâmico de interesse.
Com o intuito de um melhor aprendizado aliada à análise experimental, conta-se também com a análise via simuladores computacionais. Entretanto, simulações computacionais estão sujeitas a erros, mesmo que os dados de entrada sejam  representáveis, o resultado de uma simples operação matemática  pode não ser representável.
Em vez de um resultado verdadeiro, o computador retorna uma aproximação. Pequenos erros de arredondamento se acumulam e são propagados em sucessivos cálculos. 
Para problemas mais complexos, esses arredondamentos podem resultar em respostas totalmente erradas \cite{Galias2013,Ove2001}. 

Além de arredondamento, existem outras fontes de erros computacionais, como erro de truncamento que se origina ao truncar sequências infinitas de operações aritméticas em um número finito de etapas. No entanto, as soluções geradas por computadores são muitas vezes aceitas como respostas verdadeiras  \cite{Galias2013,Kulisch1983}. 

Simulações computacionais  também muitas vezes não coincidem com as respostas obtidas por experimentos, ficando difícil evidenciar qual resposta apresenta um resultado certo. Logo, tanto com a presença dos erros numéricos quanto dos erros intrínsecos existentes nos sistemas, a utilização da   Análise Intervalar \cite{Moore1979,Rothwell2012} é um método poderoso no controle destes erros. Métodos intervalares foram desenvolvidos com o objetivo de controlar erros de arredondamento em cálculos de ponto flutuante e seu uso cresce devido à motivação de se controlar estes erros \cite{Ruetsch2005,Moore1979}. A ideia é que em vez de usar um único valor de ponto flutuante para representar um número, o que implicaria em um erro se o número não é representável na máquina, o valor é representado por limites superior e inferior que são representáveis na máquina. Atualmente a Análise Intervalar é aplicada em várias áreas, dentre sistemas elétricos de potência \cite{Pereira2012}, processamento de sinais \cite{santana2012}, controle \cite{Banerjee2015,Yunlong2015}, entre outras.

O que se pretende neste trabalho é aplicar a Análise Intervalar a um clássico circuito elétrico RLC série, aplicando-se um degrau unitário na entrada.
Nesta análise levam-se em consideração os erros de arredondamento do computador e a tolerância dos componentes. Em seguida, pretende-se comparar os resultados com a simulação tradicional e com os dados obtidos no laboratório.

O restante do trabalho está organizado da seguinte forma. Na seção $\ref{sec:cp}$ do trabalho, são levantados conceitos básicos do artigo. Na seção $\ref{sec:met}$, a metodologia proposta é apresentada. Em seguida, na seção 4, os resultados obtidos são descritos e discutidos. Finalmente, na seção 5 são apresentadas as conclusões.

 \section{Conceitos Preliminares}
 \label{sec:cp}
\subsection{Circuito Estudado}

\begin{figure}[htb] 
		\begin{center}
			\includegraphics[angle=0, scale=1]{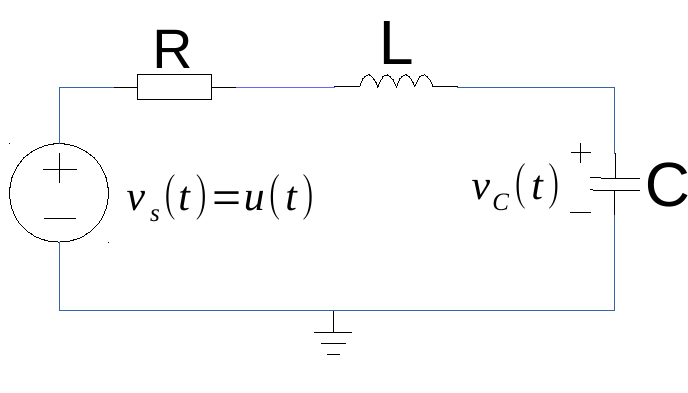} 
			\caption{Circuito RLC série.} 
			\label{fig:rlcserie}
		\end{center}
	\end{figure}

Um dos circuitos mais estudados no curso de Engenharia Elétrica é o circuito RLC série mostrado na Figura \ref{fig:rlcserie}. A sua utilização se deve ao fato de ser um circuito clássico e amplamente visto. As análises para este circuito foram realizadas de acordo com  \cite{Ogata,Sadiku}.
Inicialmente aplicando a lei de Kirchhoff das tensões no circuito, obtêm-se  a equação (\ref{eq:lkt}), consequentemente  a corrente do circuito é dada por (\ref{eq:it}),
\begin{equation}
v_{s}(t) = Ri(t) + L \frac{d}{dt}i(t) + v_{C}(t)
\label{eq:lkt}
\end{equation}
\begin{equation}
i(t) = C\frac{d}{dt}v_{C}(t).
\label{eq:it}
\end{equation}

Substituindo (\ref{eq:it}) em (\ref{eq:lkt}) têm-se que 
\begin{equation}
v_{s}(t) = RC \frac{d}{dt}v_{C}(t) + LC \frac{d^{2}}{dt^{2}}v_{C}(t) + v_{C}(t).
\label{eq:vst}
\end{equation}

Considerando os parâmetros que caracterizam uma resposta de segunda ordem e identificando-os em (\ref{eq:vst}), obtêm-se (\ref{eq:vt})
\begin{equation}
\omega_{o}^2v_{s}(t) =  \frac{d^{2}}{dt^{2}}v_{C}(t) + 2\xi\omega_{o}\frac{d}{dt}v_{C}(t) + \omega_{o}^2v_{C}(t)
\label{eq:vt}
\end{equation}
onde  $\xi$ é o fator de amortecimento ou frequência de neper, expressa em nepers
por segundo [Np/s], sendo $\xi = \frac{R}{2}\sqrt{\frac{C}{L}}$, $\omega_{o}$  a
frequência natural não amortecida do circuito expressa em radianos por segundos [rad/s] dada por  $\omega_{o} = \frac{1}{\sqrt{LC}}$ e  $\omega_d = \omega_o\sqrt{1 - \xi^2}$ é conhecida como frequência natural amortecida.

Com uma entrada em degrau unitário sendo aplicada ao circuito, a tensão de saída através do capacitor no domínio da frequência é dada por (\ref{eq:ent})
\begin{equation}
V_C(s) = \frac{1}{s}\frac{\omega_o^2}{s^2 + 2\xi\omega_o s + \omega_o^2}.
\label{eq:ent}
\end{equation}
 
Escolhidos os componentes de forma que o sistema apresente uma resposta subamortecida  $(0 < \xi <1)$ por meio da transformada inversa de Laplace de (\ref{eq:ent}), obtêm-se a equação (\ref{eq:1}), que é a resposta no domínio tempo ao se aplicar um degrau unitário na entrada para o circuito da Figura \ref{fig:rlcserie}
\begin{equation}
\small v_C(t) = 1 - e^{-\xi \omega_ t}\left ( cos\omega_d t + \frac{\xi}{\sqrt{1- \xi^2}}sen\omega_d t \right).
\label{eq:1}
\end{equation}

Em circuitos de segunda ordem  com resposta subamortecida,  algumas características do transitório são analisadas para a identificação da dinâmica do sistema, sendo estas apresentadas a seguir.
\begin{itemize}
    \item  Sobressinal Máximo ($M_p$): diferença entre o valor de pico e o valor final. Pode ser calculado por:
\begin{equation}
M_p = e^\frac{-\pi\xi}{\sqrt{1-\xi^2}}100\%.
\label{eq:sob2}
\end{equation}

\item Tempo de subida ($t_s$): tempo necessário para a resposta alcançar pela primeira vez o valor de regime.

\item Tempo de pico ($t_p$): tempo necessário para a resposta alcançar o primeiro pico de sobressinal, dado por
\begin{equation}
t_p = \frac{\pi}{\omega_d}.
\label{eq:tpico}
\end{equation}

\item Tempo de acomodação ($t_A$): tempo necessário para a resposta chegar em uma faixa do valor final, em geral de $\pm 2$\% a $\pm 5\%$.

Para uma tolerância de  $\pm$ 2\%, seu valor também pode ser definido por quatro vezes a constante de tempo, como mostrado em (\ref{eq:tA}).
\begin{equation}
t_A = \frac{4}{\xi\omega_o}.
\label{eq:tA}
\end{equation}
\end{itemize}

\subsection{Análise Intervalar}

Na análise por intervalos, os números são representados por um limite inferior e um limite superior, obtendo assim um intervalo \cite{Ruetsch2005}. Intervalos são comumente denotados por letras maiúsculas, tal como X. As extremidades de  X são denotadas frequentemente por \underline{X} e $\overline{X}$, respectivamente, de modo que X = [\underline{X}, $\overline{X}$]. Se suas extremidades são iguais \underline{X} = $\overline{X}$, esse novo número é um número real \cite{Rothwell2012}. 
A interseção  de dois intervalos $X \cap Y$ é um conjunto de números reais que pertence a ambos X e Y. A união $X \cup Y$ é um conjunto de números reais que pertence a X ou Y (ou ambos). Se $X \cap Y$ não é vazio, então $X \cap Y$ e $X \cup Y$ são intervalos que podem ser calculados por
\begin{eqnarray}
X \cap Y  = [\max(\underline{X}, \underline{Y}), \min(\overline{X}, \overline{Y})]\\
X \cup Y  = [\min(\underline{X}, \underline{Y}), \max(\overline{X}, \overline{Y})].
\end{eqnarray}  

As operações intervalares de adição, subtração e multiplicação são definidas como:
\begin{eqnarray}
X + Y = [\underline{X} + \underline{Y}, \overline{X} +  \overline{Y}]\\
X - Y = [\underline{X} - \overline{Y}, \overline{X} - \underline{Y}]\\
X \cdot Y = [\min\textit{(S)}, \max\textit{(S)}]
\end{eqnarray}
onde \textit{S} é o conjunto \{${\underline{X}\underline{Y},\underline{X}\overline{Y}, \overline{X}\underline{Y}, \overline{X}\overline{Y}}$\}. Se Y não contém o número zero, então o quociente X/Y é dado por
\begin{equation}
X/Y = X \cdot (1/Y) 
\end{equation}
onde $1/Y = [1/\overline{Y},1/\underline{Y}]$.

É importante mencionar que a adição e a multiplicação são associativas e comutativas, porém a distributividade não é geral para todos os casos, ou seja, duas expressões, que são equivalentes em aritmética verdadeira podem não ser equivalentes em aritmética intervalar \cite{Nepomuceno2016}.
Por exemplo, seja $X$ um intervalo,  a função $$ f(X)=X(1-X)$$ e uma extensão intervalar de $f$  $$ F_{1}(X)=X-X\cdot X .$$ Tem-se que $ f([0,1])=[0,1].(1-[0,1])=[0,1] $ enquanto $ F_{1}([0,1])=[0,1]-[0,1].[0,1]=[-1,1] $.

\subsection{Toolbox Intlab}

Computacionalmente é utilizado o Intlab para a análise intervalar, este é um toolbox para Matlab que suporta intervalos reais e complexos, vetores e matrizes.
Ele é projetado para ser muito rápido. De fato, não é muito mais lento do que os algoritmos puros de ponto flutuante \cite{Ru99a}.

A análise intervalar implica em computação com conjuntos. O Intlab foi projetado para obter rigorosas soluções. Uma solução para um determinado problema é produzida sob a forma de um intervalo que contêm a verdadeira solução. A aritmética finita da máquina é resolvida através de arredondamentos para o ponto flutuante mais próximo, isto é, arredondando a extremidade da esquerda para o número inferior da máquina mais próximo ou igual ao ponto final exato de um intervalo, e o ponto final direito para o número mais próximo da máquina maior do que ou igual ao exato ponto final direito \cite{Rothwell2012}.

O Intlab permite operações básicas a serem realizadas em intervalos reais e complexos, escalares, vetores e matrizes. Estas operações são inseridas semelhante a aritmética real e complexa no Matlab. Por exemplo, se a matriz $ \mathbf{A} $ é introduzida, em seguida, $ A^2 $ realiza $\mathbf{A}\times \mathbf{A}$ em aritmética intervalar, enquanto que $ A.^2 $ resulta em cada componente de $\mathbf{A}$ elevado ao quadrado, usando aritmética intervalar \cite{Hargreaves02intervalanalysis}.

No Intlab a variável é vista como um intervalo. Por exemplo, seja um resistor de 100 $\Omega$ com uma variação de $\pm 5\%$, este será representado no toolbox Intlab como
\begin{verbatim}
R = intval [95;105]
\end{verbatim}
ou seja, a variável é representada de maneira que englobe o valor nominal mais a tolerância especificada do componente.

As funções padrões, tais como as funções trigonométricas e exponenciais, estão disponíveis e
são utilizadas na forma Matlab habitual.

\section{Metodologia}
\label{sec:met}

Inicialmente foi realizada a simulação tradicional do circuito da Figura \ref{fig:rlcserie} por meio do Matlab. Para obter a curva de resposta subamortecida como proposto, foram utilizados um indutor com 0,1 H e 7,8 $\Omega$, um capacitor = 100 nF, um resistor de 100 $\Omega$. A resistência do indutor foi considerada devido aos aspectos construtivos e utilizamos o valor nominal indicado no próprio indutor. Assim o valor total da resistência série foi de  $R = 107,8 \Omega$. Em seguida,  os resultados obtidos a partir dessa simulação foram comparados com a resposta para o mesmo circuito no laboratório. A Tabela \ref{table:cop} traz a lista dos componentes utilizados na execução prática do circuito.
\begin{table}[!ht]
	\caption{Lista de componentes utilizados no circuito.}
	\label{table:cop}
	\centering
	\begin{tabular}{c c }	
		\hline
		Elemento  & Tipo/Valor nominal \\ \hline 
	    Osciloscópio  & DSO-X, 202A, 70MHZ \\
	    Matriz de Contatos & EPB0058 \\
	    Resistor & 100 $\Omega$ \\
	    Indutor & 0,1 H e 7,8 $\Omega$  \\
	    Capacitor & 100 nF \\ \hline
	\end{tabular}
\end{table}

Partindo da resposta transitória obtida pelo osciloscópio foram estimados os seguintes parâmetros $t_s$, $M_p$, $t_p$. Consequentemente foram calculados os valores dos parâmetros que caracterizam a dinâmica do circuito:   $\xi$,  $\omega_{o}$  e  $\omega_d$. Em seguida esses valores foram comparados com os previamente simulados.
Com os resultados encontrados foi feita a simulação por meio do toolbox Intlab onde os componentes do circuito são vistos como intervalos, tal que $R = [R - \delta_R, R + \delta_R]$, $L = [L - \delta_L, L + \delta_L]$, $C = [C - \delta_C, C + \delta_C]$, em que $\delta$ corresponde às tolerâncias relacionadas a cada componentes. Para o indutor é considerado uma tolerância de  $\pm 10\%$, para o capacitor $\pm 20\%$, para o resistor e para a resistência interna do indutor  $\pm 5\%$.  



 \section{Resultados}
 
\begin{figure}[!htb]
	\centering
	\includegraphics[width=8cm,height=6cm]{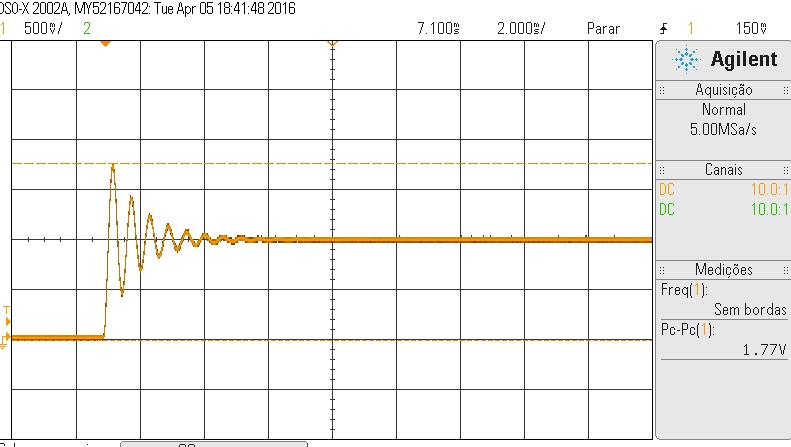}
	\caption{Resposta ao degrau obtida experimentalmente pelo osciloscópio.}
	\label{fig:osc}
\end{figure}

A Figura \ref{fig:osc} mostra a resposta, obtida por meio de um osciloscópio, ao aplicar um degrau unitário na entrada de um circuito RLC série (Figura \ref{fig:rlcserie}).
 
\begin{figure}[htb]
	\centering
	\includegraphics[width=8.5cm,height=6cm]{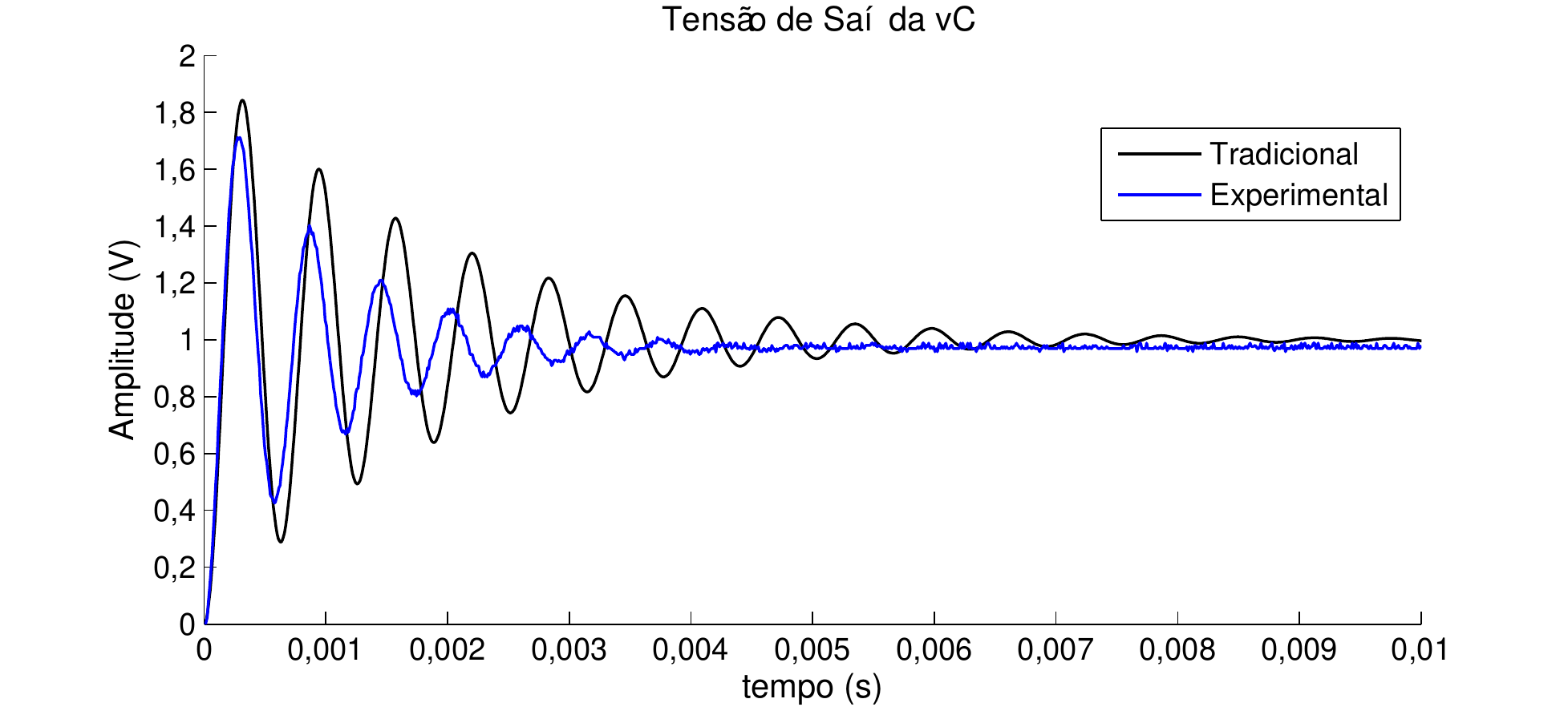}
	\caption{Respostas ao degrau obtidas experimentalmente e por simulação.}
	\label{fig:exata}
\end{figure}

A Figura \ref{fig:exata} mostra a comparação entre a simulação da curva tradicional referente a equação (\ref{eq:sob2}) e a resposta obtida experimentalmente, mostrada na Figura \ref{fig:osc}.

A partir da Figura \ref{fig:exata} é possível perceber que ambas as curvas também não coincidem totalmente e apresentam um erro aparentemente considerável. O sistema experimental apresentou um maior coeficiente de amortecimento, consequentemente diminuindo a frequência amortecida $\omega_d$. Considerando os erros existentes associados às propriedades inerentes dos componentes, a precisão finita dos equipamentos de medição e os erros propagados durante a simulação, é difícil afirmar qual representa fielmente a característica do circuito. Sendo assim é mais confiável utilizar uma simulação que contenha os erros associados.

\begin{figure}[htb]
	\centering
	\includegraphics[width=8.5cm,height=6cm]{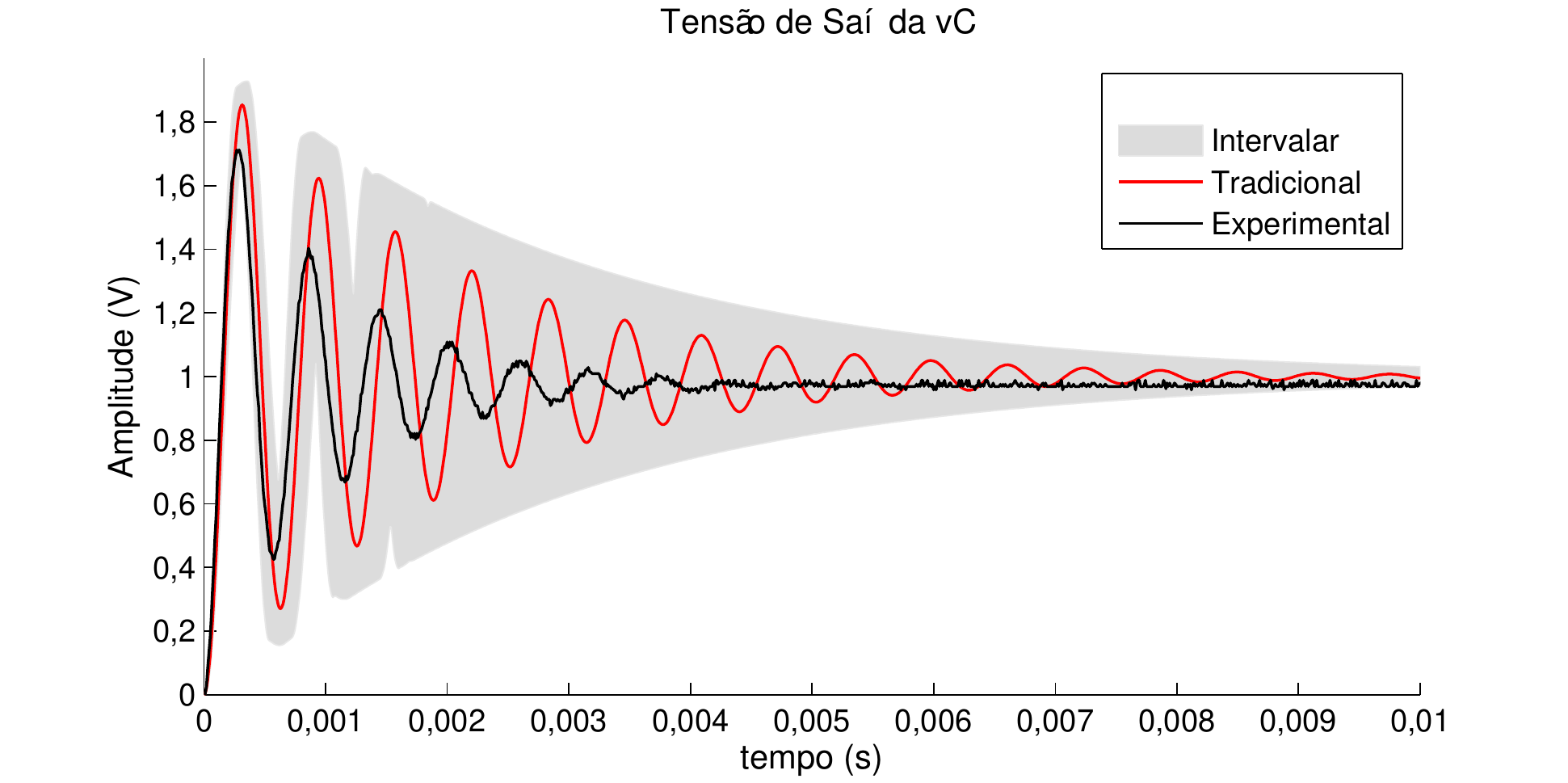}
	\caption{Comparação entre as respostas ao  degrau intervalar, por simulação tradicional e experimental.}
	\label{fig:int}
\end{figure}

Utilizando o toolbox Intlab foram  obtidas as respostas apresentadas na Figura \ref{fig:int}. É possível perceber que a resposta encontrada através do toolbox Intlab engloba  a simulação tradicional e a resposta obtida experimentalmente. Analisando as características da resposta da  curva obtida em laboratório e comparando com as simulações da resposta tradicional e da resposta intervalar, tem-se a Tabela \ref{tab:esp},  que contêm as especificações da resposta transitória obtida pela simulação tradicional, pela  resposta experimental e pela resposta intervalar.
	\begin{table}[htb]
	\small
	\caption{Especificações da resposta transitória.}
	\label{tab:esp}
	\centering
	\begin{tabular}{c c c c}
		\hline
	         & Simulação  & Experimento  & Intlab \\  \hline 
	  $ M_p$ & 0,8440    & 0,7125    & [0,6656; 0,9230]\\
	  $ t_s$ & 0,00016 & 0,000146 & [0,000137; 0,00019]\\
	  $ t_p $& 0,000314   & 0,000315 & [0,000267; 0,00036] \\
	  \hline
    \end{tabular}
\end{table}	

A partir das especificações da resposta transitória mostradas na Tabela \ref{tab:esp}, foram calculados os parâmetros que caracterizam a dinâmica do sistema. Indicados na Tabela \ref{tab:comp}.
\begin{table}[htb]
	\small
	\caption{Parâmetros que caracterizam a dinâmica do sistema.}
	\label{tab:comp}
	\centering
	\begin{tabular}{c c c c}
		\hline
	
	  & Simulação  & Experimento  & Intlab \\ \hline
	  $\xi$ & 0,0539 &  0,10729 & [0,02548; 0,12848]\\ 
	  $\omega_d$ & 9985,5 & 9951,196 & [8685,2; 11773,9]\\
	  $\omega_o$& 10000 & 10008,97 & [8703,8; 11785,1]\\	\hline
	\end{tabular}
\end{table}	
	
Ao analisar as  Tabelas \ref{tab:esp} e \ref{tab:comp}, é possível perceber que  os parâmetros que caracterizam a dinâmica do sistema, bem como os parâmetros que caracterizam a resposta transitória estão dentro dos intervalos estimados pelo Intlab. E, a partir da Figura \ref{fig:int} percebe-se que as curvas da resposta experimental e da resposta da simulação tradicional estão inclusas na resposta obtida pelo Intlab, como era esperado. 

\section{Conclusões}
 
O presente artigo apresenta uma metodologia para o tratamento de incertezas na simulação de sistemas dinâmicos, com foco em circuitos elétricos, em que as variáveis de entrada e de saída são modeladas por intervalos. Uma vez que respostas por simulações tradicionais e por experimentos em laboratórios apresentam diferenças consideráveis devido aos erros associados a simulação e aos erros inerentes do circuito prático, como pode ser visto na Figura \ref{fig:exata}. 

Com isso a solução da simulação computacional apresentada é útil para a análise de sistemas, pois esta encontra intervalos que contêm as respostas experimentais bem como as respostas obtidas via simulações tradicionais, garantindo os resultados. Caracterizando a possibilidade de se trabalhar as respostas obtidas em laboratórios, com simulações de forma que os métodos abordados durante o ensino se tornem mais condizentes, satisfazendo as respostas assim encontradas durante os estudos.

A análise intervalar desenvolvida neste trabalho se mostra eficaz e apresenta resultados satisfatórios. Na mesma linha afirmada por \cite{Rothwell2012}, a incorporação de incertezas por intervalos pode-se constituir em um método simples e eficiente para apresentar simulações tecnicamente e didaticamente mais coerentes com os dados experimentais. Acredita-se que a análise intervalar possa ser um conteúdo curricular importante para cursos de Engenharia.

\section*{Agradecimentos}
Agradecemos à CAPES, CNPq/INERGE, FAPEMIG e à Universidade Federal de São João del-Rei pelo apoio.


\bibliography{exemplo,revisao}
\end{document}